\documentclass[12pt]{amsart}
\usepackage{}

\usepackage{amsmath,amssymb}
\usepackage{amsfonts}
\usepackage{amsthm}
\usepackage{latexsym}
\usepackage{graphicx}
\usepackage{epsfig}

\newtheorem{theorem}{Theorem}[section]

\newtheorem{corollary}[theorem]{Corollary}

\newtheorem{prop}{Proposition}[section]
\newtheorem{rema}[prop]{Remark}

\makeatletter \@addtoreset{equation}{section} \makeatother

\def\({\left(}
\def\){\right)}

\begin{document}

\title{Interpolating between constrained Li-Yau and Chow-Hamilton
Harnack inequalities for a nonlinear parabolic equation}

\author{Jia-Yong Wu $^*$ }


\address{Department of Mathematics, Shanghai Maritime University,
Haigang Avenue 1550, Shanghai 201306, P. R. China}
\email{jywu81@yahoo.com}

\renewcommand{\subjclassname}{%
  \textup{2000} Mathematics Subject Classification}

\date{\today}


\maketitle

\markboth{Jia-Yong Wu}
{Interpolating between constrained Li-Yau
and Chow-Hamilton Harnack inequalities}

\begin{abstract}
We establish a one-parameter family of Harnack inequalities
connecting the constrained trace Li-Yau differential Harnack inequality
for a nonlinear parabolic equation to the constrained
trace Chow-Hamilton Harnack inequality for
this nonlinear equation with respect to evolving metrics related to
Ricci flow on a 2-dimensional closed manifold. This result
can be regarded as a nonlinear version of the previous work of
Y. Zheng and the author (Arch. Math. 94 (2010), 591-600).
\end{abstract}

\section{\textbf{Introduction}}

Let $(M^2,g(t))$, $t\in [0, T)$, be a solution to the $\varepsilon$-Ricci flow on
a 2-dimensional closed manifold $M^2$. In this paper, we will establish
an interpolation between the constrained trace Li-Yau differential
Harnack inequality for a nonlinear parabolic equation with respect
to static metrics and the constrained trace Chow-Hamilton Harnack
inequality for the nonlinear parabolic equation with respect to
evolving metrics related to Ricci flow.
More precisely, given any nonnegative constant $\varepsilon$, we say that
$g(t)$ is a solution to the $\varepsilon$-{\em Ricci flow}
on a surface $M^2$ if
\begin{equation}\label{Ricfleq}
\frac{\partial}{\partial t}g_{ij}=-\varepsilon R\cdot g_{ij},
\end{equation}
where $R$ is the scalar curvature of $g(t)$. When $\varepsilon=1$,
the $\varepsilon$-Ricci flow becomes the Ricci flow. Along the $\varepsilon$-Ricci
flow, we have
\begin{equation}\label{evolu1}
\frac{\partial R}{\partial t}=\varepsilon(\Delta R+R^2).
\end{equation}
Using the maximum principle, one can see that $R\geq c$ for some
$c\in \mathbb{R}$ is preserved along the $\varepsilon$-Ricci flow.
Under the $\varepsilon$-Ricci flow, in this paper we shall study the
Harnack inequalities for the following nonlinear parabolic equation
\begin{equation}\label{heat}
\frac{\partial f}{\partial t}=\Delta f-f\ln f+\varepsilon Rf,
\end{equation}
where $\Delta$ is the Laplacian, evolved by the $\varepsilon$-Ricci flow.
Using the maximum principle, one can see that the solutions to the
nonlinear equation \eqref{heat} will
remain positive along the $\varepsilon$-Ricci flow.

\vspace{.1in}

The motivation to study the nonlinear parabolic equation \eqref{heat}
under the $\varepsilon$-Ricci flow comes from the study of expanding
Ricci solitons, which has been nicely explained in \cite{[CaoZhang]}.
The gradient expanding Ricci solitons of the Ricci flow, which arise 
from the singularity analysis of the Ricci flow, are defined to be 
complete manifolds $(M,g)$ that the following equation
\begin{equation}\label{soliton}
R_{ij}+\nabla_i\nabla_jw=cg_{ij},
\end{equation}
holds for some Ricci soliton potential $w$ and negative constant $c$.
Taking the trace of both sides of \eqref{soliton} yields
\begin{equation}\label{soliTr}
R+\Delta w=cn.
\end{equation}
Using the contracted Bianchi identity,
\begin{equation}\label{soliBich}
R-2cw+|\nabla w|^2=\rm{constant}.
\end{equation}
From \eqref{soliTr} and \eqref{soliBich}, by
choosing a proper constant, we conclude that
\begin{equation}\label{comb}
|\nabla w|^2=\Delta w-|\nabla w|^2-R+4cw.
\end{equation}

Recall that the Ricci flow solution for a complete gradient Ricci soliton
(see Theorem 4.1 in \cite{[CLN]}) is
the pull back by $\varphi(t)$ of $g$ up to the scale factor $c(t)$:
\[
g(t)=c(t)\cdot\varphi(t)^{*}g,
\]
where $c(t):=-2ct+1>0$ and $\varphi(t)$ is the 1-parameter family of
diffeomorphisms generated by $\frac{1}{c(t)}\nabla_g w$.
Then the corresponding Ricci soliton potential $\varphi(t)^{*}w$ satisfies
\[
\frac{\partial}{\partial t} \varphi(t)^{*}w=\left|\nabla\varphi(t)^{*}w\right|^2.
\]
Note that along the Ricci flow, \eqref{comb} becomes
\[
|\nabla\varphi(t)^{*}w|^2=\Delta\varphi(t)^{*}w-|\nabla\varphi(t)^{*}w|^2
-R+\frac{4c}{c(t)}\cdot\varphi(t)^{*}w.
\]
Hence the Ricci soliton potential $\varphi(t)^{*}w$ satisfies
the evolution equation
\[
\frac{\partial\varphi(t)^{*}w}{\partial t}
=\Delta\varphi(t)^{*}w-|\nabla\varphi(t)^{*}w|^2
-R+\frac{4c}{c(t)}\cdot\varphi(t)^{*}w.
\]
If we let $\varphi(t)^{*}w=-\ln \tilde{f}$, this equation becomes
\begin{equation}\label{comb3}
\frac{\partial\tilde{f}}{\partial t}=\Delta\tilde{f}+R\tilde{f}
+\frac{4c}{c(t)}\cdot \tilde{f}\ln \tilde{f}.
\end{equation}
Notice that \eqref{comb3} and \eqref{heat} with $\varepsilon=1$ are closely related
and only differ by their last terms.

Indeed, some years ago, L. Ma~\cite{[Ma]}
proved local gradient estimates for positive solutions to the
elliptic equation
\begin{equation}\label{ellform}
\Delta f-af\ln f-bf=0,
\end{equation}
where $a$ and $b$ are real constants, on a complete manifold with 
respect to static metrics. Again, we point out that equation 
\eqref{ellform} is also related to Ricci solitons. In fact, using 
\eqref{soliTr} and \eqref{soliBich}, we deduce that
\[
\Delta w-|\nabla w|^2+2cw=\rm{constant}.
\]
So equation \eqref{ellform} can be achieved by letting $w=-\ln f$.
Later Y. Yang~\cite{[Yang]} derived local gradient estimates
for positive solutions to the corresponding nonlinear
parabolic equation
\begin{equation}\label{form}
\frac{\partial f}{\partial t}=\Delta f-af\ln f-bf
\end{equation}
on a static complete manifold (see also \cite{[ChCh]}, 
\cite{[HuangMa]}, \cite{[Wu1]}, \cite{[Wu2]}). Recently, 
Yang's result has been generalized by L. Ma~\cite{[Ma2], [Ma3]}. 
We also note that in~\cite{[Hsu]}, S.-Y. Hsu proved local gradient 
estimates for the nonlinear parabolic equation \eqref{form} with 
respect to evolving metrics related to Ricci flow. Her result is 
very similar to the Yang's gradient estimates~\cite{[Yang]}
for the static metric case. In~\cite{[CaoZhang]}, X. Cao and 
Z. Zhang derived a differential Harnack inequality for equation 
\eqref{heat} under the Ricci flow on any dimensional 
Riemannian manifold. When the dimension of the manifold 
is two, the author \cite{[Wu3]} improved their result.

\vspace{.1in}

It is well known that the study of differential Harnack inequalities originated with
the work of P. Li and S.-T. Yau~\cite{[Li-Yau]} for positive
solutions of heat equations. From then on, their Harnack inequalities
are often called Li-Yau differential Harnack inequalities. More importantly,
Li-Yau techniques were then employed by R. Hamilton, who
proved Harnack inequalities for geometric evolution equations,
especially the case of the Ricci flow~\cite{[Ham2]}. At present, there
are a large number of Harnack inequalities for various evolution equations
and their applications. The interested reader can consult the book
~\cite{[CLN]} and the recent survey \cite{ni08}.

\vspace{.1in}

On the other hand, differential Harnack inequalities for (backward) heat
equations coupled with the Ricci flow have become an important
object, which were first studied by R. Hamilton \cite{[Ham1]}.
One of the excellent important work is that G. Perelman \cite{[Perelman]}
derived differential Harnack inequalities for the fundamental solution
to the conjugate heat equation coupled with the Ricci flow without
any curvature assumption. Later X. Cao~\cite{[Caox]}, and S.-L. Kuang
and Qi S. Zhang~\cite{[KuZh]} both extended Perelman's result to the
case of all \emph{positive} solutions to the conjugate heat equation
under the Ricci flow on closed manifolds with nonnegative scalar curvature.
Besides the above work, there were also many research papers (see for
example~\cite{[BaCaoPu]}, ~\cite{[CaoxHa]},~\cite{[ChTaYu]},~\cite{CCG3},
~\cite{[Guen]},~\cite{[Liu]},~\cite{[WuZheng]} and~\cite{[Zhang]}).

\vspace{.1in}

In order to make a clear statement of our Harnack inequalities,
we need to recall some known results,
which are more or less related to our results. In~\cite{[ChHam]},
B. Chow and R. Hamilton extended Li-Yau differential Harnack inequality
~\cite{[Li-Yau]} for the heat equation on a closed
manifold, which they called a constrained trace Harnack inequality.\\

\noindent \textbf{Theorem A} (Chow-Hamilton~\cite{[ChHam]}). \emph{
Let $M^n$ be a closed manifold with nonnegative Ricci
curvature. If $S$ and $T$ are two solutions to the heat equations
\[
\frac{\partial S}{\partial t}=\Delta S\quad\text{and}\quad
\frac{\partial T}{\partial t}=\Delta T
\]
with $|T|<S$, then
\[
\frac{\partial}{\partial t}\ln S-|\nabla\ln S|^2+\frac{n}{2t}
=\Delta\ln S+\frac{n}{2t}>\frac{|\nabla h|^2}{1-h^2},
\]
where $h:=T/S$.}

\vspace{.1in}

Furthermore they generalized Hamilton's trace Harnack inequality~\cite{[Ham1]}
for the Ricci flow on surfaces with positive
scalar curvature, and proved the following
constrained linear trace Harnack inequality.

\vspace{.1in}

\noindent \textbf{Theorem B} (Chow-Hamilton~\cite{[ChHam]}).
\emph{Let $g(t)$ be a solution to the Ricci flow on a closed surface $M^2$
with scalar curvature $R>0$. If $S$ and $T$ are two solutions to
\[
\frac{\partial S}{\partial t}=\Delta S+RS\quad\text{and}\quad
\frac{\partial T}{\partial t}=\Delta T+RT
\]
with $|T|<S$, then
\[
\frac{\partial }{\partial t}\ln S-|\nabla \ln S|^2+\frac 1t=\Delta
\ln S+R+\frac 1t>\frac{|\nabla h|^2}{1-h^2},
\]
where $h:=T/S$.}

\vspace{.1in}

Recently, Y. Zheng and the author~\cite{[WuZheng]} generalized Theorem B
and Chow's interpolated Harnack inequality~\cite{[Chow3]} and proved
the interpolated and constrained linear trace Harnack inequality.

\vspace{.1in}

\noindent \textbf{Theorem C} (Wu-Zheng~\cite{[WuZheng]}). \emph{
Let $g(t)$ be a solution to the $\varepsilon$-Ricci
flow~\eqref{Ricfleq} on a closed surface $M^2$ with $R>0$. If~$S$
and $T$ are solutions to the following equations
\[
\frac{\partial S}{\partial t}=\Delta S+\varepsilon RS\quad\text{and}\quad
\frac{\partial T}{\partial t}=\Delta T+\varepsilon RT
\]
with $|T|<S$, then
\[
\frac{\partial }{\partial t}\ln S-|\nabla \ln S|^2+\frac 1t=\Delta
\ln S+\varepsilon R+\frac 1t>\frac{|\nabla h|^2}{1-h^2},
\]
where $h:=T/S$.}

\vspace{.1in}

In Theorem C, if we let $T\equiv0$ , then theorem recovers the Chow's
interpolated Harnack inequality~\cite{[Chow3]}. Chow's interpolation
trick was also adapted to proving a matrix Li-Yau-Hamilton
estimate for K$\ddot{a}$hler-Ricci flow in the work of Ni~\cite{ni07}.

Very recently, the author~\cite{[Wu3]} derived an interesting interpolated
Harnack inequality for the nonlinear parabolic equation
\eqref{heat}, also extending Chow's interpolated Harnack inequality.

\vspace{.1in}

\noindent \textbf{Theorem D} (Wu~\cite{[Wu3]}).
\emph{Let $(M,g(t))$, $t\in[0,\kappa)$, be a solution to the
$\varepsilon$-Ricci flow \eqref{Ricfleq} on a closed surface with
$R>0$. Let $f$ be a positive solution to the nonlinear parabolic
equation \eqref{heat}. Then for all time $t\in(0,\kappa)$,
\[
\frac{\partial}{\partial t}\ln f-|\nabla\ln f|^2+\ln f+\frac
1t=\Delta\ln f+\varepsilon R+\frac 1t\geq 0.
\]}

\vspace{.1in}

The main purpose of this paper is to generalize Theorems C and D, and establish
an interpolated phenomenon for the nonlinear parabolic
equation \eqref{heat} under the $\varepsilon$-Ricci flow. We will see that
this interpolated Harnack inequality is very similar to that of Theorem C.
The main difference is that the parabolic equation of this paper possesses
the additional nonlinear term: $f\ln f$. Hence in this case, the proof
is a little subtle. Let $S$ and $T$ be solutions to
the following nonlinear parabolic equations
\begin{equation}\label{linear1}
\frac{\partial S}{\partial t}=\Delta S-S\ln S+\varepsilon RS
\end{equation}
and
\begin{equation}\label{linear2}
\frac{\partial T}{\partial t}=\Delta T-T\ln T+\varepsilon RT,
\end{equation}
respectively, where $\Delta$ is the Laplacian of the metric moving under the
$\varepsilon$-Ricci flow, with the property that initially
\[
0<c_0S<T<S,
\]
where $c_0$ is a free parameter, satisfying $0<c_0<1$. Note that the
above inequality is preserved along the $\varepsilon$-Ricci flow.
In fact using \eqref{linear1} and \eqref{linear2}, we compute the
evolution equation of $h=T/S$:
\begin{equation}\label{evohfu1}
\frac{\partial h}{\partial t}=\Delta h+2\nabla h\cdot\nabla\ln S-h\ln h.
\end{equation}
Applying the maximum principle to this equation, one can prove that the inequality:
$c_0<h<1$ (and hence $c_0S<T<S$) is preserved under the $\varepsilon$-Ricci flow.

\vspace{.1in}

Now we give the following interpolation theorem.
\begin{theorem}\label{main}
Let $g(t)$ be a solution to the $\varepsilon$-Ricci
flow~\eqref{Ricfleq} on a closed surface $M^2$ with the initial
scalar curvature satisfying
\begin{equation}\label{condition}
R(g(0))\geq -\frac{2\ln c_0}{1-c_0^2}-1>0,
\end{equation}
where $c_0$ is a free parameter, satisfying $0<c_0<1$. If~$S$ and $T$ are solutions to~\eqref{linear1} and~\eqref{linear2} with $0<c_0S<T<S$ (this condition preserved by the $\varepsilon$-Ricci flow), then
\begin{equation}\label{intercons1}
\frac{\partial }{\partial t}\ln S-|\nabla \ln S|^2+\ln S+\frac 1t=\Delta
\ln S+\varepsilon R+\frac 1t>\frac{|\nabla h|^2}{1-h^2},
\end{equation}
where $h:=T/S$.
\end{theorem}

\begin{rema}
We would like to compare with Theorem C above. Theorem
\ref{main} can be regarded as a nonlinear version of Theorem C.
In Theorem \ref{main}, if we remove the nonlinear terms: $S\ln S$ in
\eqref{linear1} and $T\ln T$ in \eqref{linear2}, then the term:
$\ln S$ in \eqref{intercons1} will disappear, and we can immediately get
Theorem C under a slightly stronger scalar curvature assumption.
\end{rema}

\begin{rema}
The theorem is also true on complete noncompact surface
when the maximum principle can be used. For example, we can
assume that the solution to the $\varepsilon$-Ricci flow
is complete with the curvature and all the covariant
derivatives being uniformly bounded, and $\Delta\ln S$
has a lower bound for all time $t$.
\end{rema}

As a consequence of Theorem \ref{main}, we have a classical Harnack inequality.
\begin{theorem}\label{main2}
Let $g(t)$, $t\in(0,\kappa)$ be a solution to the
$\varepsilon$-Ricci flow~\eqref{Ricfleq} on a closed surface $M^2$
with the initial scalar curvature satisfying \eqref{condition}.
Let $S$ and $T$ be two solutions to~\eqref{linear1} and~\eqref{linear2}
with $0<c_0S<T<S$.  Assume that $(x_1,t_1)$ and
$(x_2,t_2)$, $0<t_1<t_2$, are two points in $M^2\times(0,\kappa)$. Let
\[
\Gamma:=\frac 14\inf_{\gamma}\int^{t_2}_{t_1}
e^t\left(\left|\frac{d\gamma}{dt}(t)\right|^2+\frac 4t\right)dt,
\]
where $\gamma$ is any space-time path joining $(x_1,t_1)$ and
$(x_2,t_2)$. Then we have
\[
e^{t_1}\ln S(x_1,t_1)<e^{t_2}\ln S(x_2,t_2)+\Gamma.
\]
\end{theorem}

The rest of this paper is organized as follows. In
Section~\ref{sec2}, we will prove Theorem~\ref{main}. The proof
nearly follows the proof of \cite{[WuZheng]}, which needs a
lengthy but straight-forward computation and makes use of the
parabolic maximum principle. In Section~\ref{sec3}, using
Theorem~\ref{main}, we will prove Theorem~\ref{main2}
by the standard arguments.

\section{\textbf{Proof of Theorem \ref{main}}}\label{sec2}

Under the $\varepsilon$-Ricci flow~\eqref{Ricfleq}, we can compute that
\begin{equation}\label{evolubian}
\frac{\partial }{\partial t}\ln
S=\Delta \ln S+|\nabla \ln S|^2-\ln S+\varepsilon R
\end{equation}
and
\begin{equation}\label{evolu2}
\frac{\partial}{\partial t}(\Delta)=\varepsilon R\Delta,
\end{equation}
where the Laplacian $\Delta$ is acting on smooth functions. Now we
can finish the proof of Theorem \ref{main}.
\begin{proof}[Proof of Theorem \ref{main}]
The proof follows from a direct computation and the parabolic
maximum principle. Here we mainly follow the arguments of \cite{[WuZheng]}.
Note that the equation \eqref{heat} is nonlinear.
So our case is a little more complicated. Let
\[
Q:=\Delta \ln S+\varepsilon R=\frac{\partial }{\partial t}\ln
S-|\nabla \ln S|^2+\ln S,
\]
where $S$ is a positive solution to the equation~\eqref{linear1}.
Following \cite{[WuZheng]}, using \eqref{evolubian} and \eqref{evolu2}
we compute that
\begin{equation*}
\begin{aligned}
\frac{\partial Q}{\partial t}&=\Delta\left(\frac{\partial}{\partial
t}\ln S\right)+\left(\frac{\partial}{\partial t}\Delta\right)\ln
S+\varepsilon\frac{\partial R}{\partial t}\\
&=\Delta\left(\Delta\ln S+|\nabla\ln S|^2-\ln S+\varepsilon
R\right)+\varepsilon R\Delta\ln S+\varepsilon\frac{\partial
R}{\partial t}\\
&=\Delta Q+\Delta|\nabla\ln S|^2+(\varepsilon R-1)Q+\varepsilon R-\varepsilon^2
R^2+\varepsilon\frac{\partial R}{\partial t},
\end{aligned}
\end{equation*}
where we used the equations~\eqref{evolu1}, \eqref{evolubian} and
~\eqref{evolu2}. Using the Bochner formula,
\begin{equation*}
\begin{aligned}
\frac{\partial Q}{\partial t}
&=\Delta Q+2|\nabla\nabla\ln S|^2+2\nabla\Delta\ln S\cdot\nabla\ln S
+R|\nabla\ln S|^2\\
&\quad+(\varepsilon R-1)Q+\varepsilon R-\varepsilon^2R^2
+\varepsilon\frac{\partial R}{\partial t}\\
&=\Delta Q+2|\nabla\nabla\ln S|^2
+2\nabla Q\cdot\nabla\ln S+R|\nabla\ln S|^2\\
&\quad-2\varepsilon\nabla R\cdot\nabla\ln S
+(\varepsilon R-1)Q+\varepsilon R-\varepsilon^2R^2
+\varepsilon\frac{\partial R}{\partial t}\\
&=\Delta Q+2\nabla Q\cdot\nabla\ln S-(\varepsilon R+1)Q
+2\left|\nabla\nabla\ln S+\frac{\varepsilon}{2}Rg\right|^2\\
&\quad+R\left|\nabla\ln S-\varepsilon\nabla\ln R\right|^2
+\varepsilon R[\varepsilon(\Delta\ln R+R)]+\varepsilon R.
\end{aligned}
\end{equation*}
Hence
\begin{equation}
\begin{aligned}\label{evolQ1}
\frac{\partial Q}{\partial t}&\geq\Delta Q+2\nabla Q\cdot\nabla\ln
S-(\varepsilon R+1)Q+2\left|\nabla\nabla\ln
S+\frac{\varepsilon}{2}Rg\right|^2\\
&\quad+\varepsilon R[\varepsilon(\Delta\ln R+R)]
+\varepsilon R.
\end{aligned}
\end{equation}

Next by \eqref{evohfu1}, the evolution equation of $\nabla h$ is given by
\begin{equation}
\begin{aligned}\label{evograhfu1}
\frac{\partial}{\partial t}(\nabla h)&=\nabla \left(\frac{\partial
h}{\partial t}\right)\\
&=\nabla\left[\Delta h+2\nabla h\cdot\nabla\ln S-h\ln h\right]\\
&=\Delta\nabla h+2\langle\nabla\nabla\ln S,\nabla h\rangle
+2\langle\nabla\ln S,\nabla\nabla h\rangle-\frac{R\nabla h}{2}-(1+\ln h)\nabla h.
\end{aligned}
\end{equation}
Under the $\varepsilon$-Ricci flow, using~\eqref{evograhfu1}, we have
\begin{equation*}
\begin{aligned}
\frac{\partial}{\partial t}|\nabla h|^2 &=2\nabla
h\left(\frac{\partial}{\partial t}\nabla h\right)
-g^{ki}g^{lj}\frac{\partial}{\partial t}g_{kl}\nabla_i h\nabla_j h\\
&=2\nabla h\left[\Delta\nabla h+2\langle\nabla\nabla\ln S,\nabla
h\rangle +2\langle\nabla\ln S,\nabla\nabla h\rangle-\frac{R\nabla
h}{2}-(1{+}\ln h)\nabla h\right]\\
&\quad+\varepsilon R|\nabla h|^2\\
&=\Delta|\nabla h|^2-2|\nabla\nabla h|^2+4\langle\nabla\nabla\ln
S,\nabla h\nabla h\rangle+2\langle\nabla\ln S,\nabla|\nabla
h|^2\rangle\\
&\quad+\left[(\varepsilon-1)R-2(1+\ln h)\right]|\nabla h|^2.
\end{aligned}
\end{equation*}
We also compute
\[
\frac{\partial}{\partial t}(1-h^2)=\Delta(1-h^2)+2\langle\nabla\ln
S,\nabla(1-h^2)\rangle+2|\nabla h|^2+2h^2\ln h.
\]

Next we shall compute the evolution equation of \texttt{$\frac{|\nabla h|^2}{1-h^2}$}.
Recall the following general result that if two
functions $E$ and $F$ satisfy the heat equations of the form
\[
\frac{\partial E}{\partial t}=\Delta E+A\quad \mathrm{and}
\quad
\frac{\partial F}{\partial t}=\Delta F+B,
\]
where $A$ and $B$ are some functions, then
\[
\frac{\partial}{\partial t}\left(\frac EF\right)=\Delta \left(\frac
EF\right)+\frac{2}{F^2}\langle\nabla E,\nabla
F\rangle-\frac{2E}{F^3}|\nabla F|^2+\frac AF-\frac{EB}{F^2}.
\]
Applying this result to
\[
E:=|\nabla h|^2,\quad F:=1-h^2,
\]
\[
B:=2\langle\nabla\ln
S,\nabla(1-h^2)\rangle+2|\nabla h|^2+2h^2\ln h
\]
and
\begin{equation*}
\begin{aligned}
A:&=-2|\nabla\nabla h|^2+4\langle\nabla\nabla\ln S,\nabla h\nabla
h\rangle+2\langle\nabla\ln S,\nabla|\nabla
h|^2\rangle\\
&\quad+\left[(\varepsilon-1)R-2(1+\ln h)\right]|\nabla h|^2,
\end{aligned}
\end{equation*}
we get that
\begin{equation*}
\begin{aligned}
\frac{\partial}{\partial t}\left(\frac{|\nabla h|^2}{1-h^2}\right)
&=\Delta\left(\frac{|\nabla h|^2}{1-h^2}\right)
+\frac{2\langle\nabla(1-h^2),\nabla|\nabla h|^2\rangle}{(1-h^2)^2}
-\frac{2|\nabla h|^2}{(1-h^2)^3}\left|\nabla(1-h^2)\right|^2\\
&\quad+\frac{1}{1-h^2}\cdot\left[-2|\nabla\nabla
h|^2+4\langle\nabla\nabla\ln S,\nabla h\nabla h\rangle\right]\\
&\quad+\frac{2}{1-h^2}\cdot\langle\nabla\ln S,\nabla|\nabla
h|^2\rangle+\frac{(\varepsilon-1)R-2(1+\ln h)}{1-h^2}|\nabla h|^2\\
&\quad-\frac{2|\nabla h|^2}{(1-h^2)^2}\cdot
\left[\langle\nabla\ln S,\nabla(1-h^2)\rangle+|\nabla h|^2+h^2\ln h\right].
\end{aligned}
\end{equation*}
Rearranging terms yields
\begin{equation}
\begin{aligned}\label{evograhfu}
\frac{\partial}{\partial t}\left(\frac{|\nabla h|^2}{1-h^2}\right)
&=\Delta\left(\frac{|\nabla h|^2}{1-h^2}\right)
+2\left\langle\nabla\left(\frac{|\nabla
h|^2}{1-h^2}\right),\nabla\ln S\right\rangle\\
&\quad-\frac{2}{(1-h^2)^3}\left|2h\nabla h\nabla
h+(1-h^2)\nabla\nabla h\right|^2\\
&\quad+\frac{4}{1-h^2}\langle\nabla\nabla\ln S,\nabla h\nabla
h\rangle-\frac{2|\nabla h|^4}{(1-h^2)^2}\\
&+\frac{(\varepsilon-1)R-2(1+\ln h)}{1-h^2}|\nabla h|^2
-\frac{2h^2\ln h}{(1-h^2)^2}|\nabla h|^2.
\end{aligned}
\end{equation}

Thus we define
\[
P:=Q-\frac{|\nabla h|^2}{1-h^2}=\Delta \ln S+\varepsilon
R-\frac{|\nabla h|^2}{1-h^2}.
\]
Combining~\eqref{evolQ1} and~\eqref{evograhfu}, we conclude that
\begin{equation*}
\begin{aligned}
\frac{\partial}{\partial t}P&\geq\Delta P+2\nabla P\cdot\nabla\ln
S-(\varepsilon R+1)Q+\varepsilon R+2\left|\nabla\nabla\ln
S+\frac{\varepsilon}{2}Rg\right|^2\\
&\quad+\varepsilon R[\varepsilon(\Delta\ln
R+R)]+\frac{2}{(1-h^2)^3}\left|2h\nabla h\nabla
h+(1-h^2)\nabla\nabla h\right|^2\\
&\quad-\frac{4}{1-h^2}\langle\nabla\nabla\ln S,\nabla h\nabla
h\rangle+\frac{2|\nabla h|^4}{(1-h^2)^2}\\
&\quad+\frac{(1-\varepsilon)R+2(1+\ln h)}{1-h^2}|\nabla h|^2
+\frac{2h^2\ln h}{(1-h^2)^2}|\nabla h|^2\\
&=\Delta P+2\nabla P\cdot\nabla\ln S+2\left|\nabla\nabla\ln S
+\frac{\varepsilon}{2}Rg
-\frac{\nabla h\nabla h}{1-h^2}\right|^2\\
&\quad+\varepsilon R[\varepsilon(\Delta\ln R+R)]
+\frac{2}{(1-h^2)^3}\left|2h\nabla h\nabla h +(1-h^2)\nabla\nabla
h\right|^2\\
&\quad-(\varepsilon R+1)Q+\varepsilon R
+\frac{(1+\varepsilon)R+2(1+\ln h)}{1-h^2}|\nabla h|^2
+\frac{2h^2\ln h}{(1-h^2)^2}|\nabla h|^2.
\end{aligned}
\end{equation*}
Hence we have
\begin{equation}
\begin{aligned}\label{guanjian}
\frac{\partial}{\partial t}P
&\geq\Delta P+2\nabla P\cdot\nabla\ln S+P^2-(\varepsilon
R+1)P+\varepsilon R[\varepsilon(\Delta\ln R+R)]\\
&\quad+\varepsilon R
+\frac{|\nabla h|^2}{1-h^2}\left(R+1+\frac{2\ln h}{1-h^2}\right),
\end{aligned}
\end{equation}
where we used the elementary inequality
\[
\left|\nabla\nabla\ln S+\frac{\varepsilon}{2}Rg-\frac{\nabla h\nabla
h}{1-h^2}\right|^2\geq\frac 12\left(\Delta\ln S+\varepsilon
R-\frac{|\nabla h|^2}{1-h^2}\right)^2=\frac {P^2}{2}.
\]
Since $0<c_0<h<1$ and the function
$\frac{2\ln h}{1-h^2}$ is increasing on $(0,1)$, then
\[
\frac{2\ln h}{1-h^2}>\frac{2\ln c_0}{1-c_0^2}.
\]
By the assumption of the theorem, using the maximum principle,
we can see that the inequality \eqref{condition} still
holds under the $\varepsilon$-Ricci flow. Hence
\[
R+1+\frac{2\ln h}{1-h^2}>R+1+\frac{2\ln c_0}{1-c_0^2}>0
\]
for all time $t$. Therefore, since $(M,g(t))$ has positive
scalar curvature, \eqref{guanjian} becomes
\[
\frac{\partial}{\partial t}P\geq\Delta P+2\nabla P\cdot\nabla\ln S
+P^2-(\varepsilon R+1)P+\varepsilon R[\varepsilon(\Delta\ln R+R)].
\]
Adding $\frac 1t$ to $P$ yields
\begin{equation*}
\begin{aligned}
\frac{\partial}{\partial t}\left(P+\frac 1t\right)&\geq\Delta
\left(P+\frac 1t\right)+2\nabla\left(P+\frac 1t\right)\cdot\nabla\ln
S+\left(P+\frac 1t\right)\left(P-\frac 1t\right)\\
&\quad-(\varepsilon R+1)\left(P+\frac 1t\right)+\varepsilon
R\left[\varepsilon(\Delta\ln R+R)+\frac 1t\right].
\end{aligned}
\end{equation*}
Recall that the trace Harnack inequality for the
$\varepsilon$-Ricci flow on a closed surface
proved by B. Chow in \cite{[Chow3]} (see also Lemma 2.1 in
\cite{[WuZheng]}) implies
\[
\frac{\partial\ln R}{\partial t}-\varepsilon|\nabla\ln
R|^2=\varepsilon(\Delta\ln R+R)\geq-\frac 1t,
\]
since $g(t)$ has positive scalar curvature. Hence
\begin{equation*}
\begin{aligned}
\frac{\partial}{\partial t}\left(P+\frac 1t\right)&\geq\Delta
\left(P+\frac 1t\right)+2\nabla\left(P+\frac 1t\right)\cdot\nabla\ln
S+\left(P+\frac 1t\right)\left(P-\frac 1t\right)\\
&\quad-(\varepsilon R+1)\left(P+\frac 1t\right).
\end{aligned}
\end{equation*}
It is clear to see that
\[
P+1/t>0.
\]
for very small positive $t$. Then applying the maximum principle
to the above evolution formula, we conclude that
\[
P+1/t>0
\]
for all positive time $t$, and hence the desired theorem follows.
\end{proof}

For Theorem \ref{main}, if we let $\varepsilon=0$, then
\begin{corollary}\label{corogd}
Let $M^2$ be a closed surface with the scalar curvature satisfying
\eqref{condition}. If $S$ and $T$ are solutions to
\[
\frac{\partial S}{\partial t}=\Delta S-S\ln S
\quad and\quad
\frac{\partial T}{\partial t}=\Delta T-T\ln T
\]
with $0<c_0S<T<S$, then
\[
\frac{\partial }{\partial t}\ln S-|\nabla \ln S|^2+\ln S+\frac 1t=\Delta
\ln S+\frac 1t>\frac{|\nabla h|^2}{1-h^2},
\]
where $h:=T/S$.
\end{corollary}

\vspace{.1in}

If we set
\[
\bar{g}=\varepsilon^{-1}g \quad\text{and}\quad
\alpha=\varepsilon^{-1}
\]
in Theorem \ref{main}, then
\[
\bar{\Delta}=\varepsilon\Delta\quad\text{and}\quad
\bar{R}=\varepsilon R.
\]
Hence Theorem \ref{main} can be rephrased as follows:
\begin{corollary}\label{coro}
Let $\bar{g}(t)$ be a solution to the Ricci flow
on a closed surface $M^2$ with the initial
scalar curvature satisfying
\[
\alpha \bar{R}(\bar{g}(0))\geq-\frac{2\ln c_0}{1-c_0^2}-1>0,
\]
where $\alpha$ is a positive constant and $c_0$ is a free parameter,
satisfying $0<c_0<1$. If $S$ and $T$ are solutions to
\[
\frac{\partial S}{\partial t}=\alpha\bar{\Delta} S-S\ln S+\bar{R}S
\quad and \quad
\frac{\partial T}{\partial t}=\alpha\bar{\Delta} T-T\ln T+\bar{R}T
\]
with $0<c_0S<T<S$, then
\[
\frac{\partial }{\partial t}\ln S-\alpha|\bar{\nabla}\ln S|^2+\ln S+\frac 1t=\alpha\bar{\Delta}
\ln S+\bar{R}+\frac 1t>\frac{\alpha|\bar{\nabla} h|^2}{1-h^2},
\]
where $h:=T/S$.
\end{corollary}

\section{Proof of Theorem~\ref{main2}}\label{sec3}

In the rest of this paper, we will prove Theorem~\ref{main2}
by using Theorem~\ref{main}. The proof is quite standard
by integrating the inequality \eqref{intercons1}.
We include it here for completeness.
\begin{proof}[Proof of Theorem~\ref{main2}]
We pick a space-time path $\gamma(x,t)$ joining $(x_1,t_1)$ and
$(x_2,t_2)$ with $t_2>t_1>0$. Along $\gamma$, by
Theorem \ref{main} we have
\begin{equation*}
\begin{aligned}
\frac{d}{dt}\ln S(x,t)&=\frac{\partial}{\partial t}\ln S+\nabla\ln
S\cdot\frac{d\gamma}{dt}\\
&>|\nabla \ln S|^2-\ln S-\frac 1t+\frac{|\nabla h|^2}{1-h^2}+\nabla\ln
S\cdot\frac{d\gamma}{dt}\\
&\geq-\frac14\left|\frac{d\gamma}{dt}(t)\right|^2-\ln S-\frac 1t.
\end{aligned}
\end{equation*}
Hence
\[
\frac{d}{dt}\left(e^t\ln S(x,t)\right)
>-e^t\left(\frac14\left|\frac{d\gamma}{dt}(t)\right|^2+\frac 1t\right).
\]
Integrating this inequality from the time $t_1$ to $t_2$ yields
\[
e^{t_1}\ln S(x_1,t_1)-e^{t_2}\ln S(x_2,t_2)
<\int^{t_2}_{t_1}e^t\left(\frac14\left|\frac{d\gamma}{dt}(t)\right|^2
+\frac 1t\right)dt.
\]
By the definition of $\Gamma$, we finish the proof of Theorem~\ref{main2}.
\end{proof}

\section*{Acknowledgments}
The author would like to thank an anonymous referee for many valuable
suggestions which lead to improve the paper. This work is
partially supported by NSFC (No. 11101267) and the Science and Technology
Program of Shanghai Maritime University (No. 20120061).

\end{document}